\documentclass[10pt]{article}

\usepackage{amsmath,amsfonts, latexsym,graphicx}

\newcommand{\U}{{\mathcal U}}
\newcommand{\0}{{\mathbf 0}}
\newcommand{\C}{{\mathbb C}}
\newcommand{\Z}{{\mathbb Z}}

\newcommand{\Q}{{\mathbb Q}}
\newcommand{\N}{{\mathbb N}}
\newcommand{\cL}{{\mathbb L}}

\newcommand{\W}{{\cal W}}
\newcommand{\strat}{{\mathfrak S}}
\newcommand{\Proj}{{\mathbb P}}
\newcommand{\hyp}{{\mathbb H}}
\newcommand{\supp}{\operatorname{supp}}

\newcommand{\im}{\mathop{\rm im}\nolimits}

\newcommand{\Adot}{\mathbf A^\bullet}
\newcommand{\Bdot}{\mathbf B^\bullet}
\newcommand{\Cdot}{\mathbf C^\bullet}

\newtheorem{defn0}{Definition}[section]
\newtheorem{prop0}[defn0]{Proposition}
\newtheorem{conj0}[defn0]{Conjecture}
\newtheorem{thm0}[defn0]{Theorem}
\newtheorem{lem0}[defn0]{Lemma}
\newtheorem{corollary0}[defn0]{Corollary}
\newtheorem{example0}[defn0]{Example}
\newtheorem{remark0}[defn0]{Remark}
\newtheorem{question0}[defn0]{Question}

\newenvironment{defn}{\begin{defn0}\hskip -.06in .}{\end{defn0}}
\newenvironment{prop}{\begin{prop0}\hskip -.06in .}{\end{prop0}}

\newenvironment{thm}{\begin{thm0}\hskip -.06in .}{\end{thm0}}
\newenvironment{lem}{\begin{lem0}\hskip -.06in .}{\end{lem0}}
\newenvironment{cor}{\begin{corollary0}\hskip -.06in .}{\end{corollary0}}
\newenvironment{exm}{\begin{example0}\hskip -.06in .\rm}{\end{example0}}
\newenvironment{rem}{\begin{remark0}\hskip -.06in .\rm}{\end{remark0}}

\newcommand{\defref}[1]{Definition~\ref{#1}}
\newcommand{\propref}[1]{Proposition~\ref{#1}}
\newcommand{\thmref}[1]{Theorem~\ref{#1}}
\newcommand{\lemref}[1]{Lemma~\ref{#1}}
\newcommand{\corref}[1]{Corollary~\ref{#1}}

\newcommand{\secref}[1]{Section~\ref{#1}}
\newcommand{\remref}[1]{Remark~\ref{#1}}

\newcommand{\qed}{\mbox{$\Box$}}
\newenvironment{proof}{\noindent {\bf Proof.}}{\qed\vskip 6pt}

\title{Vanishing Cycles and Thom's $a_f$ Condition\footnote{\mbox{   }   AMS subject classifications 32B15, 32C35, 32C18, 32B10.
\newline   \mbox{   } \mbox{   }   Keywords: Thom's $a_f$ condition, vanishing cycles, micro-support.}}

\author{David B. Massey}

\date{}

\begin{document}

\baselineskip = 14pt

\maketitle

\begin{abstract} We give a complete description of the relationship between the vanishing cycles of a complex of sheaves along a function $f$ and Thom's $a_f$ condition.
\end{abstract}

\sloppy

%\newpage

%\tableofcontents

%\newpage

\section{Introduction}\label{sec:intro}

Let $\U$ be an open subset of $\C^{n+1}$, and let $\tilde f:\U\rightarrow\C$ be a complex analytic function. We let $\Sigma \tilde f$ denote the critical locus of $\tilde f$. Suppose that $M$ and $N$ are complex submanifolds of $\U$.

\smallskip

{\it Thom's $a_{\tilde f}$ condition} (see, for instance, \cite{mather}) is a relative Whitney (a) condition. The $a_{\tilde f}$ is condition is important for several reasons. First, it is an hypothesis of Thom's second isotopy lemma, which allows one to conclude that maps trivialize; see \cite{mather}. Second, the $a_{\tilde f}$ condition, and the existence of stratifications in which all pairs of strata satisfy the $a_{\tilde f}$ condition, is essential in arguments such as that used by L\^e in \cite{relmono} to prove that Milnor fibrations exist even when the domain is an arbitrarily singular space. Third, the $a_{\tilde f}$ condition is closely related to constancy of the Milnor number in families of isolated hypersurface singularities; see \cite{lesaito} and below.

There are at least two important general results about the $a_{\tilde f}$ condition. There is the above-mentioned existence of $a_{\tilde f}$ stratifications, proved first in the affine setting above by Hamm and L\^e, following an argument of F. Pham,  in Theorem 1.2.1 of  \cite{hammlezariski}, and then in a different manner for an arbitrary analytic domain by Hironaka in \cite{hironakastratflat}, and there is the theorem that Whitney stratifications in which $V(\tilde f):=\tilde f^{-1}(0)$ is a union of strata  are $a_{\tilde f}$ stratifications, proved independently by Parusi\'nski  in \cite{parusw_f}, and Brian\c con, P. Maisonobe, and M. Merle in \cite{bmm}.

\smallskip

We wish to formulate the $a_{\tilde f}$ condition in conormal terms. So, we need a preliminary definition.

\begin{defn} The relative conormal space $T^*_{\tilde f_{|_M}}\U$ is given by 
$$
T^*_{\tilde f_{|_M}}\U :=\{(x, \eta)\in T^*\U\ |\ \eta(T_xM\cap \ker d_x\tilde f)=0\}.
$$
\end{defn}

\begin{rem} Note that $T^*_{\tilde f_{|_M}}\U$ equals the conormal space $T^*_M\U:=\{(x, \eta)\in T^*\U\ |\ \eta(T_xM)=0\}$ if and only if $d(\tilde f_{|_M})$ has constant rank zero, i.e.,  if  and only if $\tilde f$ is locally constant on $M$.
\end{rem}

Now, we can give the conormal definition of Thom's $a_{\tilde f}$ condition.

\begin{defn} The {\bf pair $(M, N)$ satisfies Thom's $a_{\tilde f}$ condition at a point $x\in N$} if and only if there is an inclusion, of fibers over $x$, $\big(\overline{T^*_{\tilde f_{|_M}}\U}\big)_x\subseteq \big(T^*_{\tilde f_{|_N}}\U\big)_x$.

The {\bf pair $(M, N)$ satisfies Thom's $a_{\tilde f}$ condition} if and only if it satisfies the $a_{\tilde f}$ condition at each point $x\in N$.
\end{defn}

\begin{rem} Note that if $\tilde f$ is a locally constant function, then the $a_{\tilde f}$ condition reduces to condition (a) of Whitney.
\end{rem}

\bigskip

In this paper, we prove what is essentially a generalization of the result of L\^e and Saito in \cite{lesaito}; let us recall this result, and then give the formulation which generalizes nicely.

Let $(z_0, \dots, z_n)$ be coordinates on $\U$, let $Y:=\U\cap (\C\times\{\0\})$, and assume that $Y\subseteq V(\tilde f)$.  For small $a\in\C$, define the family $\tilde f_a: (\U\cap V(z_0-a), \0)\rightarrow(\C, 0)$ by $\tilde f_a(z_1, \dots, z_n):=\tilde f(a, z_1, \dots, z_n)$. Assume that $\dim_\0\Sigma\tilde f_0 =0$.

\begin{thm}\label{thm:lesaito1}{\rm (L\^e-Saito, \cite{lesaito})}  For all small $a$, the Milnor number $\mu_{\tilde f_a}(\0)$ is independent of $a$ if and only if the only component of the critical locus of $\tilde f$, $\Sigma \tilde f$, containing the origin is $Y$ and $(\U-Y, Y)$ satisfies Thom's $a_{\tilde f}$ condition at $\0$.
\end{thm}

\begin{rem}\label{rem:nonsplit} We remark that, in the above setting, if $(\U-Y, Y)$ satisfies Thom's $a_{\tilde f}$ condition at $\0$, then the only component of  $\Sigma \tilde f$ containing the origin is $Y$. However, the proof of this requires the non-splitting result proved independently by  Gabrielov \cite{gabrielov}, Lazzeri \cite{lazzerimono}, and L\^e \cite{leacampo}.
\end{rem}

\bigskip

Using the main result of L\^e in \cite{leattach}, together with the non-splitting result of \remref{rem:nonsplit}, we can reformulate the result of L\^e and Saito as:

\begin{thm}\label{thm:lesaito2}{\rm (2nd version of L\^e-Saito Theorem, \cite{lesaito})}  For all small $a$, there is an inclusion of the Milnor fiber of $\tilde f_a$ at $(a, \0)$ into the Milnor fiber of $\tilde f$ at $\0$ which induces an isomorphism on integral cohomology if and only if  $(\U-Y, Y)$ satisfies Thom's $a_{\tilde f}$ condition at $\0$.
\end{thm}

We wish to reformulate the result of L\^e-Saito in terms of vanishing cycles. For the remainder of this paper, we let $X$ be a complex analytic subspace of $\U$, and let $f:=\tilde f_{|_X}$. 

\smallskip
Fix a base ring $R$, which is regular, Noetherian, and has finite Krull dimension, e.g., $\Z$, $\Q$, or $\C$. Let $\Adot$ be a bounded, constructible complex of sheaves of $R$-modules on $X$.

We shall use the nearby cycles, $\psi_f\Adot$, and the vanishing cycles, $\phi_f\Adot$, of $\Adot$ along $f$; we refer the reader to \cite{numcontrol}, Appendix B and \cite{dimcasheaves}. More technical references are \cite{kashschast} and \cite{schurbook}. We shall almost always include a shift by $-1$ when we apply the nearby and vanishing cycles, and we remind the reader that $\psi_f[-1]\Adot$ and $\phi_f[-1]\Adot$ are complexes of sheaves of $R$-modules on $V(f)$, with stalk cohomologies  at a point $x\in V(f)$ given by hypercohomology and relative hypercohomology of the Milnor fiber as follows:
$$
H^k(\psi_f[-1]\Adot)_x\cong \hyp^{k-1}(B_\epsilon\cap X\cap f^{-1}(a);\ \Adot),
$$
and
$$
H^k(\phi_f[-1]\Adot)_x\cong \hyp^k(B_\epsilon(x)\cap X, B_\epsilon\cap X\cap f^{-1}(a);\ \Adot),
$$
where $B_\epsilon(x)$ is a small ball (open or closed) of radius $\epsilon$ centered at $x$ in $\U$, and $0<|a|\ll\epsilon$. In the familiar case where $\Adot=\Z^\bullet_X$, this means that the stalk cohomology in degree $k$ of $\psi_f\Adot$ (respectively, $\phi_f\Z^\bullet_X$) (without the shift) at a point $x\in V(f)$ is isomorphic to the (respectively, reduced) cohomology in degree $k$ of the Milnor fiber of $f$ at $x$. We also remind the reader that, in the case where $X=\U$ and $\Adot=\Z^\bullet_{\U}$, the support of $\phi_f[-1]\Z^\bullet_{\U}$ is contained in $V(f)\cap \Sigma f$.

Below, and later, we will consider iterated vanishing cycles of the form $\phi_g[-1](\phi_f[-1]\Adot)$; when the domain of $g$ is $\U$ or $\C^{n+1}$, we shall continue to write simply $g$ in place of $g_{|_{V(f)}}$. What is the point of considering such iterated vanishing cycles?

Consider the case where $\Adot=\Z^\bullet_X$ and $g$ is the restriction to $X$ of a non-zero linear form $\mathfrak l:\C^{n+1}\rightarrow\C$. Then, the stalk cohomology $H^k(\phi_{\mathfrak l}[-1](\phi_f[-1]\Z^\bullet_X))_\0$ is isomorphic to the relative hypercohomology module
$$\hyp^k\big(B_\epsilon\cap X, B_\epsilon\cap X\cap V(\mathfrak l-a);\ \phi_f[-1]\Z^\bullet_X\big),$$
where $0<|a|\ll\epsilon\ll 1$. This module describes, on the level of cohomology, how the Milnor fibers of $f$ in a nearby hyperplane section include into the Milnor fiber at $f$ at $\0$.

\bigskip

In terms of such iterated vanishing cycles, \thmref{thm:lesaito2} becomes:

\begin{thm}\label{thm:lesaito3} {\rm (3rd version of L\^e-Saito Theorem, \cite{lesaito})}  $(\U-Y, Y)$ satisfies Thom's $a_{\tilde f}$ condition at $\0$ if and only if there exists a non-zero linear form $\mathfrak l:\C^{n+1}\rightarrow\C$ such that $H^*(\phi_{\mathfrak l}[-1]\phi_{\tilde f}[-1]\Z^\bullet_{\U})_{\0}=0$.
\end{thm}

It is this result that we generalize to the setting of arbitrary $f:X\rightarrow\C$ and with coefficients in $\Adot$. First, in \defref{def:phiconstr}, we define what it means for $\Adot$ to be {\it $\phi$-constructible} along a submanifold $M\subseteq \U$; intuitively, this notion means that the cohomology of $X$, with coefficients in $\Adot$, is ``trivial'' along $M$. Next, we define complex analytic stratifications, or merely complex analytic partitions (see \defref{def:part} and \defref{def:adotpart}), which may be weaker than Whitney stratifications, for which the Morse data associated to strata, with coefficients in $\Adot$, is still defined; a stratum in such a partition, which has a non-trivial Morse module in some degree, is called {\it $\Adot$-visible} (see \defref{def:adotpart}).

\bigskip

Our main theorem, \thmref{thm:main}, is:

\smallskip

\noindent{\bf Main Theorem}. {\it Let $M$ be a complex submanifold of $\U$ such that $M\subseteq V(f)$. Let $\W$ be a complex analytic stratification (or partition) of $X$ such that $\Adot$ is $\phi$-constructible along each stratum of $\W$.

Then, for all $\Adot$-visible $W\in\W$, $(W, M)$ satisfies the $a_f$ condition if and only if, for all $\Adot$-visible $W\in\W$, $(W, M)$ satisfies the Whitney (a) condition, and $\phi_f[-1]\Adot$ is $\phi$-constructible along $M$.}

\vskip .3in

The theorem above may look hopelessly abstract. We wish to put the reader on familiar ground by explaining what our main theorem says in the case where $\Adot$ is the constant sheaf $\Z^\bullet_X$, and the set of points where the Milnor fiber is non-trivial is $1$-dimensional (i.e., the support of the vanishing cycles is $1$-dimensional).

\smallskip

Let $M$ be a complex submanifold of $\U$ such that $M\subseteq X$. Then,  $\Z^\bullet_X$ is $\phi$-constructible along $M$ if and only if,  for all $x\in M$, for all representatives $g$ of complex analytic germs from $(X, x)$ to $(\C, 0)$ such that $x$ is a regular point of $g_{|_M}$ (i.e., such that $d_x(g_{|_M})$ is a surjection), the Milnor fiber of $g$ at $x$ has the integral cohomology of a point. In this definition, one can use simply restrictions of affine linear forms in place of more general germs $g$; we show this in \corref{cor:linver}. If $M$ is $0$-dimensional, then $\Z^\bullet_X$ is $\phi$-constructible along $M$, since the condition is vacuously satisfied. If $\dim M>0$, and $M$ is a stratum in some Whitney stratification of $X$, then $\Z^\bullet_X$ is $\phi$-constructible along $M$; however, requiring $\Z^\bullet_X$ to be $\phi$-constructible along $M$ is, in general, weaker than requiring Whitney conditions.

Now, suppose that $\W$ is a complex analytic partition of $X$ into analytic submanifolds of $\U$ (see \defref{def:part}) such that $\Z^\bullet_X$ is $\phi$-constructible along each $W\in\W$. As we may refine $\W$ to obtain a Whitney stratification, it follows that, on a generic subset of each $W\in\W$, there is a well-defined normal slice and complex link to the stratum (in the sense of Goresky and MacPherson \cite{stratmorse}). We say that a $W\in\W$ is {\it $\Z^\bullet_X$-visible} if and only if the complex link (at a generic point) of $W$ does {\bf not} have the cohomology of a point.

Let $\Sigma_{{}_{\Z}}f$ denote the {\it cohomological critical locus of $f$}, i.e., the set of points $x\in X$ such that the Milnor fiber of $f-f(x)$ at $x$ does not have the integral cohomology of a point. Suppose that $\dim V(f)\cap \overline{\Sigma_{{}_{\Z}}f} = 1$, and that $M$ is a smooth complex analytic curve contained in one of the irreducible components of $V(f)\cap \overline{\Sigma_{{}_{\Z}}f}$. Then, $\phi_f[-1]\Z^\bullet_X$ is $\phi$-constructible along $M$ if and only if, for all $x\in M$, for all representatives $g$ of complex analytic germs from $(X, x)$ to $(\C, 0)$ such that $x$ is a regular point of $g_{|_M}$, the inclusion of the Milnor fiber of $f$ at the unique point of $M\cap V(g-a)$ near $x$, for small $a\neq 0$, into the Milnor fiber of $f$ at $x$ induces an isomorphism on cohomology. In fact, in \thmref{thm:phicon}, we show that, since $\dim V(f)\cap \overline{\Sigma_{{}_{\Z}}f} = 1$, $\phi_f[-1]\Z^\bullet_X$ is $\phi$-constructible along $M$ if and only if, for all $x\in M$, there exists a single non-zero affine linear form $\mathfrak l:(\C^{n+1}, x)\rightarrow(\C, 0)$ such that $x$ is a regular point of ${\mathfrak l}_{|_M}$ and such that the inclusion of the Milnor fiber of $f$ at the unique point of $M\cap V(\mathfrak l-a)$ near $x$, for small $a\neq 0$, into the Milnor fiber of $f$ at $x$ induces an isomorphism on cohomology.

Therefore, our main theorem, \thmref{thm:main}, which we stated above, combined with \thmref{thm:phicon}, tells us, in our current situation, that the following are equivalent:

\begin{enumerate}

\item for all $\Z^\bullet_X$-visible $W\in\W$, $(W, M)$ satisfies the $a_f$ condition;
\item for all $\Z^\bullet_X$-visible $W\in\W$, $(W, M)$ satisfies the Whitney (a) condition, and for all $x\in M$, there exists a non-zero affine linear form $\mathfrak l:(\C^{n+1}, x)\rightarrow(\C, 0)$ such that $x$ is a regular point of ${\mathfrak l}_{|_M}$ and such that the inclusion of the Milnor fiber of $f$ at the unique point of $M\cap V(\mathfrak l-a)$ near $x$, for small $a\neq 0$, into the Milnor fiber of $f$ at $x$ induces an isomorphism on cohomology;
\item for all $\Z^\bullet_X$-visible $W\in\W$, $(W, M)$ satisfies the Whitney (a) condition, and for all $x\in M$, for all representatives $g$ of complex analytic germs from $(X, x)$ to $(\C, 0)$ such that $x$ is a regular point of $g_{|_M}$, the inclusion of the Milnor fiber of $f$ at the unique point of $M\cap V(g-a)$ near $x$, for small $a\neq 0$, into the Milnor fiber of $f$ at $x$ induces an isomorphism on cohomology.
\end{enumerate}
We recover the theorem of L\^e and Saito by letting $X=\U$ and $\W=\{\U\}$.

\vskip .4in

Before proving \thmref{thm:main}, we will first discuss, in \secref{sec:basic}, basic definitions and results. In \secref{sec:main}, we will prove our main theorem, and related results. Also in \secref{sec:main}, we recall results from other papers which are essential to our proofs. In \secref{sec:relations}, we shall discuss the relations between the results and techniques of this paper and those of Brian\c con, Maisonobe, and Merle in \cite{bmm}.

\section{Basic Definitions and Results}\label{sec:basic}

As in the introduction, we let $X$ be an analytic subspace of $\U$, $f:={\tilde f}_{|_X}$, and let $\Adot$ be a bounded, constructible complex of sheaves of $R$-modules on $X$. If $M$ and $N$ are complex submanifolds of $\U$, which are contained in $X$, then the $a_{\tilde f}$ condition for $(M, N)$ depends only on  $f$, and not on the extension $\tilde f$; hence, we refer simply to the $a_f$ condition.

\medskip

\begin{defn}\label{def:part} A collection $\W$ of subsets of $X$ is a {\bf (complex analytic) partition} of $X$ if and only if  $\W$ is a locally finite disjoint collection of analytic submanifolds of $\U$, which we call {\it strata}, whose union is all of $X$, and such that, for each stratum $W\in\W$, $\overline{W}$ and $\overline{W}-W$ are closed complex analytic subsets of $X$.

\smallskip

{\bf Throughout this paper, we assume that all partitions have connected strata.}

\smallskip

A partition $\W$ is a {\bf stratification} if and only if it satisfies the condition of the frontier, i.e., for all $W\in\W$, $\overline W$ is a union of elements of $\W$.
\end{defn}

\medskip

Note that, even when $\W$ is not a stratification, we nonetheless refer to elements of a partition $\W$ as strata. 

\smallskip

\begin{lem} Suppose that $\W$ and $\W^\prime$ are partitions of $X$. Let $W\in\W$. Then, there exists a unique $W^\prime\in\W^\prime$ such that $\overline{W\cap W^\prime}=\overline{W}$.
\end{lem}

\begin{proof} This is easy. Let $p\in W$. Then, by local finiteness of $\W^\prime$ and as $X=\bigcup_{W^\prime\in\W}W^\prime$, there exists an open neighborhood $\Theta$ of $p$ in $X$, and $W^\prime_1, \dots, W^\prime_d\in\W^\prime$, such that $\Theta\cap W\subseteq\bigcup_{i=1}^d(\Theta\cap W^\prime_i)$, i.e.,
$$
\Theta\cap W= \bigcup_{i=1}^d(\Theta\cap W\cap W^\prime_i).
$$

As $W$ and the $W^\prime_i$ are analytically constructible, this implies that at least one of the $\Theta\cap W\cap W^\prime_i$ is an analytically Zariski open dense subset of $\Theta\cap W$. As the elements of $\W^\prime$ are disjoint, there must be a unique such $W^\prime_i$; call it $W^\prime_p$. Now, one uses the connectedness of $W$ to conclude that $W^\prime_p$ is, in fact, the same element of $\W^\prime$ for all $p$. The desired conclusion follows.
\end{proof}

\bigskip

For most cohomological results, we do {\bf not} need a Whitney stratification of $X$ with respect to which $\Adot$ is constructible. We need merely a partition of $X$ such that the  cohomology of $X$, with coefficients in $\Adot$, is ``trivial'' along the strata. Thus, we make the following definition.

\begin{defn}\label{def:phiconstr} Let $M$ be a complex submanifold of $\U$ such that $M\subseteq X$. We say that {\bf $\Adot$ is $\phi$-constructible along $M$} if and only if,  for all $x\in M$, for all representatives $g$ of complex analytic germs from $(X, x)$ to $(\C, 0)$ such that $x$ is a regular point of $g_{|_M}$ (i.e., such that $d_x(g_{|_M})$ is a surjection), $x$ is not contained in the support, $\supp(\phi_g[-1]\Adot)$, of $\phi_g[-1]\Adot$, i.e., there exists an open neighborhood $\Theta$ of $x$ in $V(g)$ such that, for all $p\in\Theta$, $H^*(\phi_g[-1]\Adot)_p=0$.

Let $\W$ be a partition of $X$. Then, {\bf $\Adot$ is $\phi$-constructible with respect to $\W$} if and only if, for all $W\in\W$,  $\Adot$ is $\phi$-constructible along $W$.
\end{defn}

\begin{rem}\label{rem:phicon} The point of $\phi$-constructibility is that it is a purely cohomological ``replacement'' for ordinary constructibility; one which does not need to refer to a Whitney stratification. 

Of course,  if $\strat$ is a Whitney stratification of $X$, with connected strata, then  it is trivial to see that $\Adot$ is $\phi$-constructible with respect to $\strat$ if and only if $\Adot$ is constructible with respect to $\strat$.
\end{rem}

\medskip

We wish to compare $\phi$-constructibility with more standard notions. So, let $\strat$ denote a complex analytic Whitney stratification of $X$, with connected strata, with respect to which $\Adot$ is constructible. For $S\in\strat$, we let $\N_S$ and $\cL_S$ denote, respectively, the normal slice and link of the stratum $S$; see \cite{stratmorse}.

\begin{defn} A stratum $S\in\strat$ is {\bf $\Adot$-visible} if and only if the hypercohomology $\hyp^*(\N_S, \cL_S;\ \Adot)\neq 0$. We let $\strat(\Adot):= \{S\in\strat\ |\ S\ {\rm is}\ \Adot  {\text -visible}\}$.
\end{defn}
\smallskip

The point of defining $\Adot$-visible strata is that, in most cohomological results,  only the visible strata matter. In particular, if one refines $\strat$, i.e., simply throws in some extra strata, then the extra strata will be invisible; that is, the only possibly $\Adot$-visible strata in the refinement are those whose closures are equal to closures of strata in $\strat$.

\bigskip

Throughout the remainder of this paper, the {\it micro-support}, $SS(\Adot)$, of $\Adot$ will be used extensively; see \cite{kashsch}. One may also use the proposition below as the definition of $SS(\Adot)$ throughout this paper.

\bigskip

\begin{prop} {\rm (\cite{micromorse}, Theorem 4.13)} The micro-support $SS(\Adot)$ is equal to
$\displaystyle\bigcup_{S\in\strat(\Adot)}\overline{T^*_{{}_{S}}\U}$.
\end{prop}

\smallskip

Let $\tau:T^*\U\rightarrow\U$ be the projection.  For $Y\subseteq X$, we let $SS_Y(\Adot):=\tau^{-1}(Y)\cap SS(\Adot)$.

\bigskip

Now, we extend our definition of a ``visible stratum'' to certain kinds of partitions. 

\smallskip

\begin{defn}\label{def:adotpart}
A partition $\W$ of $X$ is an {\bf $\Adot$-partition} provided that 
$$
SS(\Adot)\subseteq \bigcup_{W\in\W}\overline{T^*_{{}_{W}}\U}.
$$

If $\W$ is an {\bf $\Adot$-partition}, then a stratum $W\in\W$ is {\bf $\Adot$-visible} if and only if $\overline{T^*_{{}_{W}}\U}\subseteq SS(\Adot)$. We let $\W(\Adot):= \{W\in\W\ |\ W\ {\rm is}\ \Adot{\text -visible}\}$.

Suppose that $\W$ is an {\bf $\Adot$-partition} of $X$, and $M$ is a complex submanifold of $\U$. Then,  {\bf $(\W, M)$ satisfies the $\Adot$-visible Whitney (a) condition}  (respectively, the {\bf $\Adot$-visible Thom $a_f$ condition}) if and only if, for all $\Adot$-visible $W\in\W$, $(W, M)$ satisfies Whitney's condition (a) (respectively, Thom's $a_f$ condition).

Suppose that $\W$ is an {\bf $\Adot$-partition} of $X$, and $\W^\prime$ is a partition of a closed analytic subset of $X$. Then,  {\bf $(\W, \W^\prime)$ satisfies the $\Adot$-visible Whitney (a) condition}  (respectively, the {\bf $\Adot$-visible Thom $a_f$ condition}) if and only if, for all $W^\prime\in\W^\prime$, $(\W, W^\prime)$ satisfies the $\Adot$-visible Whitney (a) condition  (respectively, the $\Adot$-visible Thom $a_f$ condition)
\end{defn}

\begin{rem} The reader should understand that the point of an $\Adot$-partition $\W$ is that, for each $\Adot$-visible stratum $S$ in $\strat$, there exists a unique $W\in\W$ such that $\overline{S}=\overline{W}$ and, hence, $\overline{T^*_{{}_{S}}\U}=\overline{T^*_{{}_{W}}\U}$. It follows at once from this, and the definition of $\Adot$-visible strata of $\W$, that, if $\W$ is an $\Adot$-partition, then
$$SS(\Adot)\ =\ \bigcup_{W\in \W(\Adot)}\overline{T^*_{{}_{W}}\U}.$$

We could, of course, define an $\Adot$-partition without using the conormal formulation in \defref{def:adotpart}. If we define the set $E(\Adot)$ of {\it $\Adot$-essential varieties} by $E(\Adot):=\{\overline{S}\ |\ S\in\strat(\Adot)\}$, then a partition $\W$ of $X$ is an $\Adot$-partition if and only if $
E(\Adot)\subseteq \{\overline{W}\ |\ W\in\W\}$. However, the conormal characterization in \defref{def:adotpart} will be very useful later.

We should also remark that in \cite{pervcohovan}, we referred to $\Adot$-partitions as {\it $\Adot$-normal partitionings}.
\end{rem}

In \cite{critpts}, we made the following definition:

\begin{defn} The $\Adot$-critical locus of $f$, $\Sigma_{{}_{\Adot}}f$, is $\{x\in X\ |\ H^*(\phi_{f-f(x)}[-1]\Adot)_x\neq 0\}$.
\end{defn}

\smallskip

The support of $\phi_f[-1]\Adot$ can be ``calculated'' as follows:

\smallskip

\begin{thm}\label{thm:vansupp} {\rm (\cite{singenrich}, Theorem 3.4)} 
$$
\supp\phi_f[-1]\Adot\ =\ V(f)\cap \overline{\Sigma_{{}_{\Adot}}f}\ =\ \{x\in V(f)\ |\ (x, d_x\tilde f)\in SS(\Adot)\}.
$$
\end{thm}

\begin{cor}\label{cor:linver} In \defref{def:phiconstr}, one may replace each reference to a complex analytic germ $g$ by the restriction to $X$ of an affine linear form and obtain a characterization of $\phi$-constructibility.

To be precise, let $M$ be a complex submanifold of $\U$ such that $M\subseteq X$. Suppose that, for all $x\in M$, for all linear forms $\mathfrak l:\C^{n+1}\rightarrow\C$  such that $x$ is a regular point of ${\mathfrak l}_{|_M}$, $x$ is not contained in $\supp(\phi_{\mathfrak l-\mathfrak l(x)}[-1]\Adot)$. Then, $\Adot$ is $\phi$-constructible along $M$.
\end{cor}

\begin{proof} By \thmref{thm:vansupp}, $x\in \supp\phi_f[-1]\Adot$ if and only if $x\in\supp\phi_{\mathfrak l-\mathfrak l(x)}[-1]\Adot$ where $\mathfrak l:= d_x\tilde f$ (where we have identified $\C^{n+1}$ with its tangent space at $x$). The corollary is immediate.
\end{proof}

\section{Main Theorems}\label{sec:main}

The following result is closely related to Proposition 8.6.4 of \cite{kashsch}.

\begin{prop}\label{prop:phicon} Let $M$ be a complex submanifold of $\U$ which is contained in $X$. Then, $\Adot$ is $\phi$-constructible along $M$ if and only if $SS_M(\Adot)\subseteq T^*_{{}_{M}}\U$.
\end{prop}
\begin{proof} By \thmref{thm:vansupp}, $p\not\in\supp\phi_g[-1]\Adot$ if and only if $d_p\tilde g\not\in SS_p(\Adot)$, where $\tilde g$ is a local extension of $g$ to $\U$. The conclusion is immediate.
\end{proof}

\smallskip

Below, we once again identify $T_p\U$ with the ambient $\C^{n+1}$, and so identify elements of $(T^*\U)_p$ with linear forms on the ambient space.

\smallskip

\begin{lem}\label{lem:genphi} Let $\W$ be a partition of $X$ such that, for all $W\in\W$ such that $\dim W>0$, for all $p\in W$, there exists a projective algebraic set $V_p\subseteq \Proj\big((T^*\U)_p\big)$ such that $\dim V_p=\dim W-1$, and such that, for all projective classes $[\mathfrak l]\in V_p$, $p$ is not contained in $\supp(\phi_{{\mathfrak l}-\mathfrak l(p)}[-1]\Adot)$.

Then, $\W$ is an $\Adot$-partition;
\end{lem}
\begin{proof} Let $S\in\strat(\Adot)$. Let $W$ be the unique element of $\W$ such that $\overline{S\cap W} = \overline{S}$. Then, $\dim W\geq\dim S$. We claim that $\dim W=\dim S$, which implies that $\overline{W}=\overline{S}$ and $\overline{T^*_{{}_{W}}\U}=\overline{T^*_{{}_{S}}\U}$; this would prove the lemma.

\smallskip

If $\dim S=n+1$, there is nothing to show. So assume that $\dim S\leq n$.

\smallskip

Suppose that $\dim W>\dim S$. Let $p\in S\cap W$.   Let $V_p$ be as in the statement of the lemma. Then, 
$$n-\dim \Proj\big((T^*_{{}_{S}}\U)_p\big) = \dim S<\dim W = \dim V_p +1,
$$
i.e., $n-1< \dim \Proj\big((T^*_{{}_{S}}\U)_p\big) + \dim V_p$. Thus, the projective algebraic subsets  $\Proj\big((T^*_{{}_{S}}\U)_p\big)$ and  $V_p$ in $\Proj\big((T^*\U)_p\big)\cong\Proj^n$ have a non-empty intersection, i.e., there exists $[\mathfrak l]\in \Proj\big((T^*_{{}_{S}}\U)_p\big) \cap V_p$. By \thmref{thm:vansupp}, $p\in\supp(\phi_{\mathfrak l-\mathfrak l(p)}[-1]\Adot)$, which contradicts that $[\mathfrak l]\in V_p$.
\end{proof}

\bigskip

\begin{defn} If $\W$ is a partition which satisfies the hypothesis of \lemref{lem:genphi}, we say that {\bf $\Adot$ is weakly $\phi$-constructible with respect to $\W$}.
\end{defn}

\begin{rem}\label{rem:weakphi} Note that \lemref{lem:genphi} enables us to talk about $\Adot$-visible strata when $\Adot$ is weakly $\phi$-constructible with respect to $\W$.
\end{rem}

\smallskip

\begin{thm}\label{thm:phicon} Let $\W$ be a partition of $X$. Then, the following are equivalent:
\begin{enumerate}
\item $\Adot$ is $\phi$-constructible with respect to $\W$;

\item $SS(\Adot)\subseteq\bigcup_{W\in\W}T^*_W\U$;

\item $\W$ is an $\Adot$-partition such that $(\W, \W)$ satisfies the $\Adot$-visible Whitney (a) condition;

\item $\Adot$ is weakly $\phi$-constructible with respect to $\W$, and $(\W, \W)$ satisfies the $\Adot$-visible Whitney (a) condition.

\end{enumerate}
\end{thm}

\begin{proof} That $(1)$ and $(2)$ are equivalent follows immediately from \propref{prop:phicon}.

\smallskip

If $\W$ is an $\Adot$-partition, then $SS(\Adot) = \bigcup_{W\in\W(\Adot)}\overline{T^*_{{}_{W}}\U}$, and $(\W, \W)$ satisfies the visible Whitney (a) condition if and only if
$$
\bigcup_{W\in\W(\Adot)}\overline{T^*_{{}_{W}}\U}\subseteq \bigcup_{W\in\W}T^*_{{}_{W}}\U.
$$
Thus, $(2)$ and $(3)$ are equivalent.

\smallskip

Now, $(1)$ and $(3)$ are equivalent, and clearly, together, they imply $(4)$. Finally, \lemref{lem:genphi} tells us that $(4)$ implies $(3)$.
\end{proof}

\bigskip

\begin{exm} In order to see why the $\Adot$-visible Whitney (a) condition is important in the above theorem, consider the following example. Let $X=\U=\C^2$, and use $y$ and $z$ for coordinates. Let $\Bdot$ be the constant sheaf (over $\Z$) on the $y$-axis, extended by zero to all of $\C^2$. Let $\Cdot$ be the constant sheaf (over $\Z$) on the $z$-axis, extended by zero to all of $\C^2$. Let $\Adot = \Bdot\oplus\Cdot$. 

Then, $SS(\Adot) = T^*_{V(z)}\U \cup T^*_{V(y)}\U$. The conormal to the origin $T^*_{\0}\U$ does {\bf not} appear in $SS(\Adot)$, because the stalk of $\Adot$ at the origin is $\Z\oplus\Z$ and so is the stalk cohomology of the Milnor fiber of a generic linear form, and the comparison map is an isomorphism; it follows that the vanishing cycles of $\Adot$ along a generic linear form are zero at the origin.

Thus, the partition $\W = \{\C^2-V(yz), V(z), V(y)-\{\0\}\}$ is an $\Adot$-partition of $\C^2$. Note that we have {\bf not} included $\{\0\}$ as a stratum. The $\Adot$-visible strata of $\W$ are $V(z)$ and $V(y)-\{\0\}$. The paragraph above tells us that $\Adot$ is weakly $\phi$-constructible with respect to $\W$. However, $(\W, \W)$ does not satisfy the $\Adot$-visible Whitney (a) condition. According to \thmref{thm:phicon}, $\Adot$ is not $\phi$-constructible with respect to $\W$.

This is easy to see in our current example. The origin is a regular point of the function given by the restriction of $y$ to $V(z)$, and yet $H^*(\phi_y\Adot)_{\0}\neq 0$, since the nearby fiber of the function $y$ is a single point which has a single $\Z$ for its cohomology (in degree $0$).
\end{exm}

\bigskip

Below, for the sake of self-containment, we state the results from other papers that we need.

\bigskip

We identify $T^*\U$ with $\U\times\C^{n+1}$. Recall that $\tau:T^*\U\rightarrow\U$ is the projection. 

\bigskip

If $\Bdot$ is a bounded, constructible complex of sheaves on $X$, and $Y\subseteq X$, we let $SS(\Bdot)_{\subseteq Y}$ denote the components of $SS(\Bdot)$ which lie above $Y$ (i.e., project by $\tau$ to a set contained in $Y$). Note that, in general, $SS(\Bdot)_{\subseteq Y}\subseteq SS_Y(\Bdot)$, and the inclusion can easily be proper. Let $i:X-V(f)\hookrightarrow X$ denote the inclusion.

The equality involving $SS(i_!i^!\Adot)_{\subseteq V(f)}$ in the next theorem can be proved using Theorem 4.2 of \cite{hypercohom}. In Proposition 4.5 of \cite{pervcohovan}, we used 4.2.1 and 3.4.2 of \cite{bmm} to prove the result when $\Adot$ is a perverse sheaf; in the perverse sheaf case, the characteristic cycle is equal to the micro-support. Our primary reason for including results about $SS(i_!i^!\Adot)_{\subseteq V(f)}$ will be discussed in the next section.

\smallskip

\begin{thm}\label{thm:fundnearby}{\rm(\cite{singenrich}, Theorem 3.2 and \cite{pervcohovan}, Proposition 4.5)} There is the following equality of subspaces of the cotangent space, $T^*\U$:
$$
SS(\psi_f[-1]\Adot) \  =\ SS(i_!i^!\Adot)_{\subseteq V(f)}\ =\  \tau^{-1}(V(f))\ \cap\ \big(\bigcup_{\substack{S\in\strat(\Adot)\\ f_{|_S}\neq{\rm\ const.}}}\overline{T^*_{f_{|_S}}\U}\big).
$$
\end{thm}

\begin{cor}\label{cor:fundnearby1} Let $\W$ be an $\Adot$-partition of $X$. Let $M$ be a complex submanifold of $\U$ such that $M\subseteq V(f)$. Then, the following are equivalent:
\begin{enumerate}
\item for all $W\in\W(\Adot)$ such that $W\not\subseteq V(f)$, $(W, M)$ satisfies the  $a_f$ condition;
\item  $\psi_f[-1]\Adot$ is $\phi$-constructible along $M$;
\item $\tau^{-1}(M)\cap SS(i_!i^!\Adot)_{\subseteq V(f)}\subseteq T^*_{{}_M}\U$.
\end{enumerate}

In addition, these equivalent conditions imply that, for all $W\in\W(\Adot)$ such that $W\not\subseteq V(f)$, $(W, M)$ satisfies the  Whitney (a) condition
\end{cor}
\begin{proof} Combine the theorem with \propref{prop:phicon} and the conormal characterization of the $a_f$ condition.
\end{proof}

\begin{cor}\label{cor:fundnearby2} Let $\W$ be an $\Adot$-partition of $X$. Let $\W^\prime$ be a Whitney (a)  partition of $V(f)$. Then, the following are equivalent:
\begin{enumerate}
\item for all $W\in\W(\Adot)$ such that $W\not\subseteq V(f)$, for all $W^\prime\in\W^\prime$, $(W, W^\prime)$ satisfies the  $a_f$ condition;
\item $\psi_f[-1]\Adot$ is $\phi$-constructible with respect to $\W^\prime$;
\item $\psi_f[-1]\Adot$ is weakly $\phi$-constructible with respect to $\W^\prime$.
\end{enumerate}
\end{cor}

\begin{proof} This follows at once from \corref{cor:fundnearby1} and  \thmref{thm:phicon}.
\end{proof}

\smallskip

\begin{rem} We wish to discuss the problem in using \corref{cor:fundnearby2} in practice; a problem that is removed by replacing the nearby cycles with the vanishing cycles.

Our primary goal in this paper, as we discussed in the introduction, is to provide a generalization of the result of L\^e and Saito, in the form given in \thmref{thm:lesaito3}. We could obtain an analogous statement, using the nearby cycles in place of the vanishing cycles, by using the equivalence of Items 1 and 3 above. The problem is that we are required to begin with a Whitney (a) partition of all of $V(f)$, instead of merely a partition of $\Sigma f$. Requiring that the smooth part of $V(f)$ satisfy the Whitney (a) condition with respect to strata of $\Sigma f$ is an unacceptable assumption, as such an assumption does not appear in the theorem of L\^e and Saito.

The way that we will fix this problem is to use the vanishing cycles, whose support is contained in the critical locus.
\end{rem}

\bigskip

We let $\pi:\U\times\C^{n+1}\times\Proj^n\rightarrow \U\times\Proj^n$ and $\nu: \U\times\Proj^n\rightarrow\U$ denote the respective projections. Recall that $\tilde f$ is our global extension of $f$ to all of $\U$ (though we could use local extensions at each point). We let $\im d\tilde f$ denote the image of $d\tilde f$ in $T^*\U$.

\begin{thm}\label{thm:blowup}{\rm(\cite{pervcohovan}, Proposition 4.3)} Suppose that $Y$ is an analytic subset of $X$. Suppose that $f$ is not constant on any irreducible component of $Y$.   Let $E$ denote the
exceptional divisor in
$\operatorname{Bl}_{\im df}\overline{T^*_{{}_{Y_{\operatorname{reg}}}}\U}\ \subseteq\ \U\times\C^{n+1}\times\Proj^n$. Suppose that $M\subseteq X$ is a complex analytic submanifold of $\U$ and that $x\in M$ is such that
$(X_{\operatorname{reg}}, M)$ satisfies Whitney's condition a) at $x$ and such that $d_x(f_{|_M})\equiv 0$.

Then, $(Y_{\operatorname{reg}}, M)$ satisfies Whitney's $a_f$ condition at $x$ if and only if there is the containment of fibres
above
$x$ given by 
$$\big(\pi(E)\big)_x\subseteq \big(\Proj(T^*_{{}_M}\U)\big)_x.$$
\end{thm}

\bigskip

Now, for each $S\in\strat$ or $S\in \W$, let $E_S$ denote the exceptional divisor of $\operatorname{Bl}_{\im d\tilde f}\overline{T^*_{{}_{S}}\U}$. Then, in our current notation, the second equality of Theorem 3.4 of \cite{singenrich} tells us:

\begin{thm}\label{thm:fundvan}{\rm(\cite{singenrich}, Theorem 3.4)} There is the following equality of subspaces of the projectivized cotangent space, $\Proj(T^*\U)$:
$$
\Proj\big(SS(\phi_f[-1]\Adot)\big) = \Proj\big(SS(\phi^\prime_f[-1]\Adot)\big) =\nu^{-1}(V(f))\ \cap\ \pi\big(\bigcup_{S\in\strat(\Adot)}E_S\big).
$$
\end{thm}

\bigskip

Note that, if $\W$ is an $\Adot$-partition, then, in \thmref{thm:fundvan}, we could have replaced $\bigcup_{S\in\strat(\Adot)}E_S$ by $\bigcup_{W\in\W(\Adot)}E_W$.

\bigskip

In previous papers, we have proved two results along the lines of our main theorem below. In Theorem 4.4 of \cite{pervcohovan}, we proved a form of this result in the case where $\Adot$ is a {\it perverse sheaf}. In the case of general $\Adot$, we proved one direction of this result in Theorem 6.5 of \cite{singenrich}. In addition to containing less general results than our current paper, both \cite{pervcohovan} and \cite{singenrich} are so abstract that the reader would have difficulty extracting the relevant results. Also, Theorem 6.5 of \cite{singenrich} is proved using Theorem 4.8 of that paper; Theorem 4.8 is misstated (though is fine in the case where it is used). For all of these reasons, we prove both directions of the theorem below.

\begin{thm}\label{thm:main} Let $\W$ be an $\Adot$-partition of $X$. Let $M$ be a complex submanifold of $\U$ such that $M\subseteq V(f)$.

Then, $(\W, M)$ satisfies the $\Adot$-visible $a_f$ condition if and only if $(\W, M)$ satisfies the $\Adot$-visible Whitney (a) condition and $\phi_f[-1]\Adot$ is $\phi$-constructible along $M$.
\end{thm}

\begin{proof}  Let $W\in\W$. From \thmref{thm:blowup}, it follows easily that: $(\dagger)$\ the pair $(W, M)$ satisfies the $a_f$ condition if and only if $(W, M)$ satisfies Whitney's condition (a) and $\nu^{-1}(M)\cap\pi(E_W)\subseteq\Proj\big(T^*_{{}_{M}}\U\big)$.

\medskip

\noindent Proof of $\Rightarrow$:

Now, suppose that $(\W, M)$ satisfies the $\Adot$-visible $a_f$ condition. Then, $(\dagger)$ immediately implies that $(\W, M)$ satisfies the $\Adot$-visible Whitney (a) condition.  Combining \thmref{thm:fundvan} with $(\dagger)$, we find that
$$
\Proj\big(SS_M(\phi_f[-1]\Adot)\big) = \nu^{-1}(M)\ \cap\ \pi\big(\bigcup_{W\in\W(\Adot)}E_W\big)\subseteq \Proj\big(T^*_{{}_{M}}\U\big).
$$ 
By \propref{prop:phicon}, this is equivalent to $\phi_f[-1]\Adot$ being $\phi$-constructible along $M$.
\bigskip

\noindent Proof of $\Leftarrow$:

Suppose that $(\W, M)$ satisfies the $\Adot$-visible Whitney (a) condition, and that $\phi_f[-1]\Adot$ is $\phi$-constructible along $M$.
Then, as above,
$$
\Proj\big(SS_M(\phi_f[-1]\Adot)\big) = \nu^{-1}(M)\ \cap\ \pi\big(\bigcup_{W\in\W(\Adot)}E_W\big)\subseteq \Proj\big(T^*_{{}_{M}}\U\big).
$$

Let $W\in\W(\Adot)$. By the $\Adot$-visible Whitney (a) condition and $(\dagger)$, what we need to show is that $\nu^{-1}(M)\cap\pi(E_W)\subseteq\Proj\big(T^*_{{}_{M}}\U\big)$, which follows from the above.
\end{proof}

\vskip .3in

Let $\phi^\prime_f[-1]\Adot$ denote the complex of sheaves on $V(f)\cap \overline{\Sigma_{{}_{\Adot}}f}$ obtained by restricting $\phi_f[-1]\Adot$ to its support. Note that $SS(\phi^\prime_f[-1]\Adot) = SS(\phi_f[-1]\Adot)$.

\begin{cor}\label{cor:phipart} Let $\W$ be an $\Adot$-partition of $X$. Let $\W^\prime$ be a Whitney (a)  partition of $V(f)\cap \overline{\Sigma_{{}_{\Adot}}f}$. Suppose that $(\W,\W^\prime)$ satisfies the $\Adot$-visible Whitney (a) condition.

Then, the following are equivalent:
\begin{enumerate}

\item $(\W, \W^\prime)$ satisfies the $\Adot$-visible $a_f$ condition;
\item  $\phi^\prime_f[-1]\Adot$ is $\phi$-constructible with respect to $\W^\prime$;
\item $\phi^\prime_f[-1]\Adot$ is weakly $\phi$-constructible with respect to $\W^\prime$;
\end{enumerate}

\noindent and, if $\W^\prime$, is, in fact, a Whitney stratification, these are equivalent to:
\begin{enumerate}
\item[4.]  $\phi^\prime_f[-1]\Adot$ is constructible with respect to $\W^\prime$.
\end{enumerate}

\end{cor}

\begin{proof} To obtain the equivalences of Items 1, 2, and 3, simply combine \thmref{thm:main} with \thmref{thm:phicon}. If $\W^\prime$, is a Whitney stratification, then the equivalence of Items 2 and 4 follows from \remref{rem:phicon}.
\end{proof}

\smallskip

\begin{exm} Let us return to the result of L\^e and Saito, which we discussed at length in the introduction. We use the assumptions and notation that we used in \thmref{thm:lesaito1}. 

Let $X=\U$, $\W=\{\U\}$, $f=\tilde f$, and $\Adot=\Z^\bullet_{\U}$. Then, the critical locus of $f$  near the origin is equal to $V(f)\cap \overline{\Sigma_{{}_{\Z^\bullet_{\U}}}f}$, which we suppose is simply $Y=\U\cap (\C\times\{\0\})$. Let $\W^\prime=\{Y\}$. 

Then, the hypotheses of \corref{cor:phipart} are satisfied. In addition, the pair $(\U, Y)$ satisfies the $a_f$ {\bf generically} along $Y$; the only question, near the origin is: ``what happens at the origin?''. As $Y$ is $1$-dimensional, to know that $\phi^\prime_f[-1]\Z^\bullet_{\U}$ is weakly $\phi$-constructible, we need to have a single non-zero linear form $\mathfrak l: \C^{n+1}\rightarrow\C$ such that  $H^*(\phi_{\mathfrak l}[-1](\phi^\prime_f[-1]\Z^\bullet_{\U}))=H^*(\phi_{\mathfrak l}[-1](\phi_f[-1]\Z^\bullet_{\U}))=0$ (recall \remref{rem:weakphi}).

Therefore, the conclusion of \corref{cor:phipart} is precisely our third version of the L\^e-Saito Theorem, which we stated in \thmref{thm:lesaito3}.
\end{exm}

\section{Relations with the Work of Brian\c con, Maisonobe, and Merle}\label{sec:relations} 

In \cite{bmm}, Brian\c con, Maisonobe, and Merle introduce the {\it condition of local, stratified triviality} -- a condition on a Whitney (a) stratification. The condition is that, for any point $x$ in a stratum $S$, every analytic transverse slice to $S$ at $x$ (of any dimension) yields a stratified homeomorphism between an open neighborhood of $x$ and the product of the slice with a open ball. See Definition 4.1.1 of \cite{bmm}.

Thus, if one has a stratification, $\W$, of $X$ of which satisfies the condition of local, stratified triviality and $\Adot$ is a bounded, constructible complex of sheaves on $X$ whose local structure depends only on the local stratified topological-type of $X$, then $\Adot$ will be $\phi$-constructible with respect to $\W$.

Consider now a Whitney stratification $\strat$ of $X$ such that $V(f)$ is a union of strata. Recall that $i:X-V(f)\hookrightarrow X$ denotes the inclusion. For each stratum $S\in \strat$ such that $S\not\subseteq V(f)$, let $\Adot_S$ denote the extension by zero, to all of $X$, of the constant sheaf $\Z^\bullet_{{}_S}$. Then, as $\strat$ is a Whitney stratification, $i_!i^!\Adot_S$ is $\phi$-constructible with respect $\strat$, and certainly $S$ is $(i_!i^!\Adot_S)$-visible.

Therefore, if $M$ is a stratum of $\strat$ and $M\subseteq V(f)$, then \propref{prop:phicon} tells us that
$$
\tau^{-1}(M)\cap SS(i_!i^!\Adot)_{\subseteq V(f)}\ \subseteq \ SS_M(i_!i^!\Adot)\ \subseteq\ T^*_{{}_M}\U,
$$
and \corref{cor:fundnearby1} tells us that the pair $(S, M)$ satisfies the $a_f$ condition.

The above is precisely the argument used in \cite{bmm} to prove that Whitney stratifications, in which $V(\tilde f):=\tilde f^{-1}(0)$ is a union of strata,  are $a_{\tilde f}$ stratifications. We remark again that this result was proved independently by Parusi\'nski  in \cite{parusw_f}. We should also remark that, because Brian\c con, Maisonobe, and Merle used characteristic cycles, instead of micro-supports, in some parts of their paper, they needed to use a perverse sheaf for our $\Adot_S$ above. Hence, rather than use the extension by zero of the constant sheaf, they used the extension by zero of the intersection cohomology complex (with constant coefficients) on $\overline{S}$.

\bigskip

The reader should understand that we included results on $\psi_f[-1]\Adot$ and $SS(i_!i^!\Adot)_{\subseteq V(f)}$ in this paper in order to show how $\phi$-constructible partitions arise in the proof of the main theorem of \cite{bmm}; most of these results appeared in some form in \cite{bmm}. However, the results of \cite{bmm} do not give us the desired generalization of the result of L\^e and Saito; for that, we need our results on the vanishing cycles in \thmref{thm:blowup}, \thmref{thm:fundvan}, \thmref{thm:main}, and \corref{cor:phipart}.

\bibliographystyle{plain}
\bibliography{Masseybib}
%\printindex
\end{document}